\documentclass{proc-l}

\usepackage{epsf}
\usepackage{epsfig}
\usepackage{graphicx}
\usepackage{xcolor}
\usepackage{bm}
\usepackage{tabularx}
\usepackage{amstext}
\usepackage{float}
\usepackage{array, nicefrac, caption}

\newcommand{\divprop}{\lfloor}
\newcommand{\1}{{\bf 1}}
\newtheorem{theorem}{Theorem}

\allowdisplaybreaks
\begin{document}

\frenchspacing

\title{Properties of the recursive divisor function and the number of ordered factorizations}
\author{T. M. A. Fink}
\address{London Institute for Mathematical Sciences, Royal Institution, 21 Albermarle St, London W1S 4BS, UK}

\date{\today}
\vspace{-0.2in}
\begin{abstract}
We recently introduced the recursive divisor function  $\kappa_x(n)$, a recursive analogue of the usual divisor function. 
Here we calculate its Dirichlet series, which is ${\zeta(s-x)}/(2 - \zeta(s))$. 
We show that $\kappa_x(n)$ is related to the ordinary divisor function by
$\kappa_x * \sigma_y = \kappa_y * \sigma_x$, where * denotes the Dirichlet convolution.
Using this, we derive several identities relating $\kappa_x$ and some standard arithmetic functions.
We also clarify the relation between $\kappa_0$ and the much-studied number of ordered factorizations $K(n)$, namely, $\kappa_0 = \1 * K$.
\vspace*{-0.25in}
\end{abstract}
\maketitle
\noindent
Several arithmetic functions, such as the divisor function and the Euler totient function, 
play a fundamental role in our understanding of the theory of numbers.
They are explicitly defined, in the sense that values of the function are not defined in terms of prior values of the function.
However, it is possible to write down meaningful recursive arithmetic functions. 
We consider two of them in this paper, and show that they are intimately related to a number of standard arithmetic functions.
\\ \indent
We recently introduced and studied the recursive divisor function \cite{Fink21,Fink22}:
	\begin{equation}
		\kappa_x(n) = n^x+ \sum_{d \divprop n} \kappa_x(d),
		\label{AH}
	\end{equation}
where $m \divprop n$ means $m \vert n$ and $m < n$.
It is the recursive analogue of the usual divisor function,
	\begin{equation*}
		\sigma_x(n) = \sum_{d|n} d^x.
	\end{equation*}
For example, 
$\kappa_0(4) = 1 + \kappa_0(1) + \kappa_0(2) = 4$ and
$\kappa_1(6) = 6 + \kappa_1(1) + \kappa_1(2) + \kappa_1(3) = 14$. 
The first 12 values of $\kappa_0$ (A067824 \cite{Sloane}) and $\kappa_1$ (A330575 \cite{Sloane}) are shown is Table 1. 
\\ \indent
While the function $\kappa_x(n)$ has received little attention \cite{Fink21,Fink22},
a related but simpler recursive arithmetic function has been studied for 90 years, namely,
the number $K(n)$ of ordered factorizations into integers greater than 1 \cite{Kalmar,Hille,Canfield,Chor,Klazara,Deleglise}.
It is defined as
	\begin{equation}
		K(n) = \varepsilon(n) + \sum_{d\divprop n} K(d),
		\label{AL}
	\end{equation} 
where $\varepsilon(n)$ is 1 for $n=1$ but zero otherwise.
For example, $K(8)=4$ because 8 is the product of integers greater than one in four ways: 
$8 = 4 \cdot 2 = 2 \cdot 4 = 2 \cdot 2 \cdot 2$.
The first 12 values of $K$ (A074206 \cite{Sloane}) are shown is Table 1.
\\ \indent
In what follows, we denote the Dirichlet convolution of two arithmetic functions $f$ and $g$ as $f * g$, 
and we denote the Dirichlet series of $f$ as $\widetilde {f}$.
The various arithmetic functions that we use in this paper are summarized in Table 1.
\section{Statement of results}
\noindent
We prove the following three theorems.
\begin{theorem}
	The Dirichlet series $\widetilde \kappa_x$ for the recursive divisor function $\kappa_x(n)$ is
	\begin{equation*}
		\widetilde \kappa_x = \sum_{n=1}^\infty \frac{\kappa_x(n)}{n^s} = \frac{\zeta(s-x)}{2 - \zeta(s)},
	\end{equation*}
	where $\zeta$ is the Riemann zeta function.
	\label{AP}
\end{theorem}
\begin{theorem}
The recursive divisor function $\kappa_x$ satisfies the following:
	\begin{align}
		\kappa_x * \sigma_y	&= \kappa_y * \sigma_x 			& & \kappa \mbox{-} \sigma \,\, {\rm exchange \,\, symmetry}					\label{BA} 	\\
		\kappa_x			&= ({\rm id}_x + \1 * \kappa_x)/2 	& & {\rm Definition \,\, of} \,\,  \kappa_x									\label{BB} \\
		\kappa_x			&= \frac{{\rm id}_x}{2} + \frac{\1*{\rm id}_x}{2^2} + \frac{\1*\1*{\rm id}_x}{2^3} + \ldots & & {\rm Series \,\, representation \,\, of} \,\,  \kappa_x \label{BC} \\
		\kappa_x			&= J_x * \kappa_0				& & {\rm Relation \,\, between} \,\,  \kappa_x	\,\, {\rm and} \,\, \kappa_0			\label{BD} \\
		\kappa_x^{-1}		&= J_x^{-1} * (2 \, \mu - \varepsilon)	& & {\rm Inverse \,\, of} \,\,  \kappa_x										\label{BE} \\
		\sigma_x			&= \kappa_x * (2 \, \1 - d)			& & {\rm Relation \,\, between} \,\,  \kappa_x	\,\, {\rm and} \,\, \sigma_x.			\label{BF} 
	\end{align}
\end{theorem}
\noindent
Note the special cases:
	$\kappa_1 = \phi * \kappa_0$;
	$\kappa_0^{-1}	= 2 \mu - \varepsilon$; and
	$\kappa_0	 = \frac{\1}{2} + \frac{\1*\1}{2^2} + \frac{\1*\1*\1}{2^3} + \ldots$.
\begin{theorem}
The number of recursive divisors is related to the number of ordered factorizations by $\kappa_0 = \1 * K$, that is, 
	\begin{align}
		\kappa_0(n) &= \sum_{d|n} K(d).
		\label{DZ}
	\end{align}
Furthermore, $K$ satisfies the following:
	\begin{align}
		K			&= (\varepsilon + \1 * K)/2									&  & {\rm Definition \,\, of} \,\,  K							\label{DA} \\
		K			&= \frac{\varepsilon}{2} + \frac{\1}{2^2} + \frac{\1*\1}{2^3} + \ldots 	&  & {\rm Series \,\, representation \,\, of} \,\, K 				\label{DB} \\
		\kappa_x		&= {\rm id}_x * K										&  & {\rm Relation \,\, between} \,\,  K 	\,\, {\rm and} \,\, \kappa_x	\label{DC} \\
		K^{-1}		&= 2 \, \varepsilon - \1									&  & {\rm Inverse \,\, of} \,\,  K.							\label{DD} 
	\end{align}
\end{theorem}
\begin{table}[b!]
\setlength{\tabcolsep}{1.6pt}
\begin{small}
\begin{tabular*}{\textwidth}{@{\extracolsep{\fill}}cllcrrrrrrrrrrrrr}
				&								&& \emph{Dirichlet}												 \\ 
				&								&& \emph{series}						&& \multicolumn{12}{c}{\emph{Values from} $n=1$ \emph{to} $n=12$} 	\\  \vspace{2pt}
$\varepsilon$		& Identity $\lfloor 1/n \rfloor$ 			&& 1									&& 1 & 0  &  0 & 0 & 0 & 0 & 0 & 0 & 0 & 0 & 0 & 0					\\  \vspace{2pt}
$\mu$			& M\"obius function					&& $ 1/\zeta(s)$						&& 1 & -1 & -1 & 0 & -1 & 1 & -1 & 0 & 0 & 1 & -1 & 0				\\  \vspace{2pt}
\1				& Constant function 1 $\equiv {\rm id}_0$	&& $\zeta(s)$							&& 1 & 1  &  1 & 1 & 1 & 1 & 1 & 1 & 1 & 1 & 1 & 1 					\\  \vspace{2pt}
${\rm id}_1$		& $1$st power function $n$			&& $\zeta(s-1)$							&& 1 & 2  &  3 & 4 & 5 & 6 & 7 & 8 & 9 & 10 & 11 & 12 				\\  \vspace{2pt}
$\phi$			& Euler totient function $\equiv J_1$		&& $\displaystyle\frac{\zeta(s-1)}{\zeta(s)}$	&& 1 & 1  &  2 & 2 & 4 & 2 & 6 & 4 & 6 & 4 & 10 & 4	 				\\  \vspace{2pt}
$d$				& No. of divisors $\equiv \sigma_0$		&& $\zeta^2(s)$						&& 1 & 2  &  2 & 3 & 2 & 4 & 2 & 4 & 3 & 4 & 2 & 6 					\\  \vspace{2pt}
$\sigma$			& Sum of divisors $\equiv \sigma_1$		&& $\zeta(s) \zeta(s-1)$					&& 1 & 3  &  4 & 7 & 6 & 12 & 8 & 15 & 13 & 18 & 12 & 28		 		\\  \vspace{2pt}
$\kappa_0$		& No. of recursive divisors 			&& $\displaystyle\frac{\zeta(s)}{2 - \zeta(s)}$	&& 1 & 2  &  2 & 4 & 2 & 6 & 2 & 8 & 4 & 6 & 2 & 16 					\\  \vspace{2pt}
$\kappa_1$		& Sum of recursive divisors 			&& $\displaystyle\frac{\zeta(s-1)}{2 - \zeta(s)}$	&& 1 & 3  &  4 & 8 & 6 & 14 & 8 & 20 & 14 & 20 & 12 & 42 			\\  \vspace{2pt}
$\kappa_0^{-1}$	& Inverse of $\kappa_0$ 				&& $\displaystyle\frac{2 - \zeta(s)}{\zeta(s)}$	&& 1 & -2 & -2 & 0 & -2 &  2 & -2 & 0 & 0 & 2 & -2 & 0 				\\  \vspace{2pt}
 $K$				& No. of ordered factorizations 			&& $\displaystyle 1/(2 - \zeta(s))$			&& 1 & 1  &  1 & 2  & 1 & 3 & 1 & 4 & 2 & 3 & 1 & 8	 				\\  \vspace{2pt}
$K^{-1}$			& Inverse of $K$ 					&& $ 2 - \zeta(s)$						&& 1 & -1 & -1 & -1 &-1 & -1 & -1 & -1 & -1 & -1 & -1 & -1 				\\  \vspace{2pt}
$\sigma_x$		& Sum of $x$th power of divisors		&& $\zeta(s) \zeta(s-x)$					&& \multicolumn{12}{l}{$1, \, 2^x+1, \, 3^x+1, \, 4^x + 2^x+1, \ldots$}		\\  \vspace{2pt}	
${\rm id}_x$ 		& $x$th power function $n^x$			&& $\zeta(s-x)$							&&  \multicolumn{12}{l}{$1, \, 2^x, \, 3^x, \, 4^x, \ldots$}				\\  \vspace{2pt}
$J_x$ 			& Jordan's totient function				&& $\displaystyle \frac{\zeta(s-x)}{\zeta(s)}$	&& \multicolumn{12}{l}{$1, \, 2^x-1, \, 3^x-1, 4^x-2^x, \ldots$}			\\  \vspace{-6pt}
$\kappa_x$ 		& Sum of $x$th power of  				&& $\displaystyle \frac{\zeta(s-x)}{2 - \zeta(s)}$	&& \multicolumn{12}{l}{$1, \, 2^x+1, \, 3^x+1, \, 4^x + 2^x+2, \ldots$}		\\  \vspace{0pt}
				& recursive divisors 					&& 							
\end{tabular*}
\end{small}
\vspace{0.1in}
\caption{\small 
For each of the arithmetic functions used in this paper, we give its Dirichlet series and the first 12 terms of its sequence
(four terms when $x$ is not specified).
}
\label{avalstable}
\end{table}
\section{Proof of Theorem 1}
\noindent
It is convenient to rewrite (\ref{AH}) as
	\begin{equation}
		2 \, \kappa_x(n) = n^x+ \sum_{d | n} \kappa_x(d).
		\label{JM}
	\end{equation}
Dividing by $n^s$ and summing over $n$,
	\begin{align*}
		2 \, \widetilde \kappa_x  = 
		2 \sum_{n=1}^\infty \frac{\kappa_x(n)}{n^s} &= \sum_{n=1}^\infty \frac{n^x}{n^s} + \sum_{n=1}^\infty \frac{1}{n^s} \sum_{d|n} \kappa_x(d).
	\end{align*}
Swapping the order of summation,	
	\begin{align*}
	2 \, \widetilde \kappa_x 
	&= \zeta(s-x) + \sum_{d=1}^\infty \kappa_x(d)  \sum_{n: \, d|n} \frac{1}{n^s} 		\\
									  &= \zeta(s-x) + \sum_{d=1}^\infty \kappa_x(d)  \sum_{n=1}^\infty \frac{1}{(d n)^s} 	\\
								 	  &= \zeta(s-x) + \sum_{d=1}^\infty \frac{\kappa_x(d)}{d^s}  \sum_{n=1}^\infty \frac{1}{n^s} 	\\
								 	  &= \zeta(s-x) + \widetilde \kappa_x  \, \zeta(s).
	\end{align*}
From this, we arrive at Theorem \ref{AP}:
	\begin{equation*}
		\widetilde \kappa_x = \frac{\zeta(s-x)}{2 - \zeta(s)}. \qed
	\end{equation*}
\section{Proof of Theorem 2}
\noindent
In what follows, we use the standard identities
$\1 * \mu = \varepsilon$,
$\1 * J_x = {\rm id}_x$ and
$\1 * {\rm id}_x = \sigma_x$.
\\ \indent
Since $\widetilde \sigma_x/\widetilde \kappa_x = \zeta(s)(2 - \zeta(s))$ is independent of $x$,  
$\widetilde \kappa_x \, \widetilde \sigma_y = \widetilde \kappa_y \, \widetilde \sigma_x$.
The Dirichlet convolutions of $\kappa_x$ and $\sigma_x$ must follow the analogous relation, so we arrive at (\ref{BA}):
	\begin{align*}
		\kappa_x * \sigma_y &= \kappa_y * \sigma_x. \qed
		\label{MD}
	\end{align*}	
We can immediately rewrite (\ref{JM}) in the form of (\ref{BB}),
	\begin{equation*}
		\kappa_x = ({\rm id}_x + 1 * \kappa_x)/2. \qed
		\label{QA}
	\end{equation*}
Iterating this recursive definition gives
	\begin{align*}
		\kappa_x &= \frac{{\rm id}_x}{2} + \frac{\1 * {\rm id}_x}{4} + \frac{\1 * \1 * \kappa_x}{4},
	\end{align*}
and so on, leading to the infinite series (\ref{BC}):
	\begin{align*}
		\kappa_x			&= \frac{{\rm id}_x}{2} + \frac{\1*{\rm id}_x}{2^2} + \frac{\1*\1*{\rm id}_x}{2^3} + \ldots. \qed
	\end{align*}
Since $\sigma_x = \1 * \1 * J_x$, from (\ref{BA}) we have $\kappa_x * J_y	= \kappa_y * J_x$.	
Setting $y=0$, and since $J_0 = \varepsilon$, we have (\ref{BD}),
	\begin{align*}
		\kappa_x &= J_x * \kappa_0. \qed
	\end{align*}
Convolving of (\ref{BB}) by $\kappa_x^{-1}$ and solving for $\kappa_x^{-1}$,
we have $\kappa_x^{-1}	=  {\rm id}_x^{-1} * (2 \varepsilon - \1)$. 
Since ${\rm id}_x = \1 * J_x$, we have ${\rm id}_x^{-1} = \mu * J_x^{-1}$, and we arrive at (\ref{BE}),
	\begin{align*}
		\kappa_x^{-1}	&= J_x^{-1} * (2 \mu - \varepsilon).  \qed
	\end{align*}
From (\ref{BE}), $J_x	 = \kappa_x * (2 \mu - \varepsilon)$. 
Convolving with $\1 * \1$, we find (\ref{BF}),
	\begin{align*}
		\sigma_x = \kappa_x * (2 \, \1 - d). \qed
	\end{align*}
\section{Proof of Theorem 3}
\noindent
Rewriting (\ref{AL}) as
	\begin{equation*}
		2 K(n) = \varepsilon(n) + \sum_{d|n} K(d),
	\end{equation*}
we immediately arrive at (\ref{DA}),
	\begin{equation*}
		K = (\varepsilon  + \1 * K)/2. \qed
	\end{equation*}
Convolving (\ref{BB}) with \1 yields $\1 * K = (\1  + \1 * \1 * K)/2$.
Replacing $1 * K$ with $\kappa_0$, we recover the definition of $\kappa_0$: $\kappa_0 = (\1  + \1 * \kappa_0)/2$.
Thus we have established (\ref{DZ}),
	\begin{equation*}
		\kappa_0 = \1 * K.  \qed
	\end{equation*}
Setting $x=0$ in (\ref{BC}), and convolving with $\mu$, we have (\ref{DB}),
	\begin{equation*}
		K = \frac{\varepsilon}{2} + \frac{\1}{2^2} + \frac{\1*\1}{2^3} + \ldots.  \qed
	\end{equation*}
Substituting $\kappa_0 = \1 * K$ into (\ref{BD}), we have (\ref{DC}),
	\begin{equation*}
		\kappa_x = {\rm id}_x * K. \qed
	\end{equation*}
Inverting $K = \mu * \kappa_0$ gives $K^{-1} = \1 * \kappa_0^{-1}$.
Inserting into this $\kappa_0^{-1} = 2 \mu - \varepsilon$ from (\ref{BE}) gives (\ref{DD}),
	\begin{equation*}
		K^{-1} = 2 \varepsilon - \1. \qed
	\end{equation*}

\begin{footnotesize}

\end{footnotesize}

\begin{thebibliography}{1}
\bibitem{Fink21} 		T. Fink, 								Recursively divisible numbers,												arxiv.org/abs/1912.07979.
\bibitem{Fink22} 		T. Fink, 								Recursively abundant and recursively perfect numbers,							arxiv.org/abs/2008.10398.
\bibitem{Sloane}		N.\ J.\ A.\ Sloane, editor, 					The On-Line Encyclopedia of Integer Sequences, published electronically at https://oeis.org, 2018.
\bibitem{Kalmar}		L. Kalmar, 							A factorisatio numerorum probelmajarol, 										\emph{Mat Fiz Lapok}		{\bf 38}, 	1 		(1931).
\bibitem{Hille} 			E. Hille,								A problem in factorisatio numerorum, 										\emph{Acta Arith} 			{\bf 2}, 	134		(1936).
\bibitem{Canfield} 		E. Canfield, P. Erd\"{o}s, C. Pomerance,		On a problem of Oppenheim concerning ``factorisatio numerorum",					\emph{J Number Theory} 		{\bf 17},  	1	 	(1983).
\bibitem{Chor} 			B. Chor, P. Lemke, Z. Mador,				On the number of ordered factorizations of natural numbers, 						\emph{Disc Math} 			{\bf 214},	123		(2000).
\bibitem{Klazara} 		M. Klazar, F. Luca,						On the maximal order of numbers in the factorisatio numerorum problem,				\emph{J Number Theory} 		{\bf 124},  470	 	(2007).
\bibitem{Deleglise}		M. Del\'{e}glise, M.\ Hernane, J.-L. Nicolas,	Grandes valeurs et nombres champions de la fonction arithm\'{e}tique de Kalm\'{a}r, 		\emph{J Number Theory} 		{\bf 128}, 	1676 	(2008).
\end{thebibliography}
\end{document}